\documentclass[12pt]{article}
\usepackage{amsfonts,amsmath,amsxtra}
\usepackage{latexsym}
\usepackage{amssymb}
\usepackage{color}
\usepackage[matrix,arrow,curve]{xy}
\def\hybrid{\topmargin 0pt      \oddsidemargin 0pt
        \headheight 0pt \headsep 0pt
        \textwidth 16.5cm
        \textheight 23cm
        \marginparwidth 0.0in
        \parskip 5pt plus 1pt   \jot = 1.5ex}
\catcode`\@=11
\def\marginnote#1{}
\newcount\hour
\newcount\minute
\newtoks\amorpm
\hour=\time\divide\hour by60 \minute=\time{\multiply\hour by60
\global\advance\minute by-\hour}
\edef\standardtime{{\ifnum\hour<12 \global\amorpm={am}%
        \else\global\amorpm={pm}\advance\hour by-12 \fi
        \ifnum\hour=0 \hour=12 \fi
      \number\hour:\ifnum\minute<10 0\fi\number\minute\the\amorpm}}
\edef\militarytime{\number\hour:\ifnum\minute<10 0\fi\number\minute}

\def\draftlabel#1{{\@bsphack\if@filesw {\let\thepage\relax
   \xdef\@gtempa{\write\@auxout{\string
      \newlabel{#1}{{\@currentlabel}{\thepage}}}}}\@gtempa
   \if@nobreak \ifvmode\nobreak\fi\fi\fi\@esphack}
        \gdef\@eqnlabel{#1}}
\def\@eqnlabel{}
\def\@vacuum{}
\def\draftmarginnote#1{\marginpar{\raggedright\scriptsize\tt#1}}

\def\draft{\oddsidemargin -0.1truein
        \def\@oddfoot{\sl WLimit.tex \hfil
        \rm\thepage\hfil\sl\today\quad\militarytime}
        \let\@evenfoot\@oddfoot \overfullrule 3pt
        \let\label=\draftlabel
        \let\marginnote=\draftmarginnote
\def\@eqnnum{{\rm (\theequation)}
\rlap{\kern\marginparsep\tt\@eqnlabel}%
\global\let\@eqnlabel\@vacuum}  }
\newfont{\Bbbb}{msbm7 scaled 1\@ptsize00}
\newcommand{\zs}{\raise-1pt\hbox{$\mbox{\Bbbb Z}$}}

scaled 1\@ptsize00
\def\numberbysection{\@addtoreset{equation}{section}
        \def\theequation{\thesection.\arabic{equation}}}
\numberbysection

\renewcommand{\theequation}{\thesection.\arabic{equation}}

\def\titlepage{\@restonecolfalse\if@twocolumn\@restonecoltrue\onecolumn
     \else \newpage \fi \thispagestyle{empty}\c@page\z@
\def\thefootnote{\fnsymbol{footnote}} }
\def\endtitlepage{\if@restonecol\twocolumn \else  \fi
        \def\thefootnote{\arabic{footnote}}
        \setcounter{footnote}{0}}  %\c@footnote\z@ }
\relax \hybrid
\makeatletter
\newdimen\normalarrayskip            % skip between lines
\newdimen\minarrayskip               % minimal skip between lines
\normalarrayskip\baselineskip \minarrayskip\jot
\newif\ifold             \oldtrue            \def\new{\oldfalse}
\def\arraymode{\ifold\relax\else\displaystyle\fi}%mode of array enrties
\def\eqnumphantom{\phantom{(\theequation)}} % ight phantom in eqnarray
\def\@arrayskip{\ifold\baselineskip\z@\lineskip\z@
     \else
     \baselineskip\minarrayskip\lineskip1\baselineskip\fi}
\def\@arrayclassz{\ifcase \@lastchclass \@acolampacol \or
\@ampacol \or \or \or \@addamp \or
   \@acolampacol \or \@firstampfalse \@acol \fi
\edef\@preamble{\@preamble
  \ifcase \@chnum
     \hfil$\relax\arraymode\@sharp$\hfil
     \or $\relax\arraymode\@sharp$\hfil
     \or \hfil$\relax\arraymode\@sharp$\fi}}
\def\@array[#1]#2{\setbox\@arstrutbox=\hbox{\vrule
     height\arraystretch \ht\strutbox
     depth\arraystretch \dp\strutbox
width\z@}\@mkpream{#2}\edef\@preamble{\halign \noexpand\@halignto
\bgroup \tabskip\z@ \@arstrut \@preamble \tabskip\z@ \cr}%
\let\@startpbox\@@startpbox \let\@endpbox\@@endpbox
  \if #1t\vtop \else \if#1b\vbox \else \vcenter \fi\fi
  \bgroup \let\par\relax
  \let\@sharp##\let\protect\relax
  \@arrayskip\@preamble}
\def\eqnarray{\stepcounter{equation}%
              \let\@currentlabel=\theequation
              \global\@eqnswtrue
              \global\@eqcnt\z@
              \tabskip\@centering              %formulae  centering
              \let\\=\@eqncr
              $$%
            \halign to \displaywidth  \bgroup
             \eqnumphantom \@eqnsel
      \hskip\@centering                               %right tab%
    $\displaystyle  \tabskip\z@ {##}$%
    &\global\@eqcnt\@ne \hskip 2\arraycolsep
         $ \displaystyle  \arraymode{##}$\hfil
    &\global\@eqcnt\tw@ \hskip 2\arraycolsep
         $\displaystyle\tabskip\z@{##}$\hfil
         \tabskip\@centering
    &{##}\tabskip\z@\cr}
\makeatother
\def\IC{\mathbb{C}}

\def\IL{\mathbb{L}}

\def\IR{\mathbb{R}}
\def\IZ{\mathbb{Z}}
\def\CA {\mathcal{A}}

\def\CC {\mathcal{C}}

\def\CE {\mathcal{E}}
\def\CF {\mathcal{F}}

\def\CH {\mathcal{H}}
\def\CI {\mathcal{I}}

\def\CL {\mathcal{L}}

\def\CS {\mathcal{S}}

\def\CU {\mathcal{U}}
\def\CV {\mathcal{V}}
\def\CW {\mathcal{W}}

\def\fb{{\frak b}}
\def\fh{{\frak h}}
\def\fH{{\frak H}}
\def\fL{{\frak L}}
\def\fn{{\frak n}}
\def\a {{\alpha}}

\def\g {{\gamma}}

\def\la{\lambda}

\def\e{\epsilon}
\def\pr {\partial}

\def\wh{\widehat}

\def\lt{\tilde{l}}
\def\Lie{{\rm Lie}}
\def\Mat{{\rm Mat}}

\def\Span{{\mathop{\rm span}}}

\def\frak{\mathfrak}

\def\gl{\mathfrak{gl}}

\def\<{\langle}
\def\>{\rangle}

\newtheorem{te}{Theorem}[section]%Usage:\begin{te}Statement\end{te}
\newtheorem{de}{Definition}[section]
\newtheorem{prop}{Proposition}[section]           %  ETC ...

\newtheorem{lem}{Lemma}[section]
\newtheorem{ex}{Example}[section]

\newcommand\rk{\operatorname{rank}}

\newcommand{\proof}{\noindent {\it Proof}. }
\newcommand\bqa{\begin{eqnarray}}
\newcommand\eqa{\end{eqnarray}}
\def\be{\begin{eqnarray}\new\begin{array}{cc}}
\def\ee{\end{array}\end{eqnarray}}
\def\beq{\begin{equation}}
\def\eeq{\end{equation}}
\def\bse{\begin{subequations}}                %%%SUBEQUATIONS
\def\ese{\end{subequations}}
\def\bp{\begin{pmatrix}}
\def\ep{\end{pmatrix}}

\newcounter{pac}[section]

\newcounter{pacc}[subsection]

 0
\begin{document}
%\draft

\title{\bf On a matrix element representation of the GKZ hypergeometric functions}
\author{A.A. Gerasimov, D.R. Lebedev and S.V. Oblezin}
\date{\today}
\maketitle

\renewcommand{\abstractname}{}

\begin{abstract}
\noindent {\bf Abstract}. We develop a representation theory
approach to the study of generalized hypergeometric functions of
Gelfand, Kapranov and Zelevisnky (GKZ). We show that the GKZ
hypergeometric functions may be identified with matrix elements of
non-reductive Lie algebras $\fL_N$ of oscillator type. The Whittaker
functions associated with principal series representations of
$\gl_{\ell+1}(\IR)$ being  special cases  of GKZ hypergeometric
functions, thus admit along with a standard  matrix element
representations associated with reductive Lie algebra
$\gl_{\ell+1}(\IR)$, another matrix element representation in terms
of $\fL_{\ell(\ell+1)}$.

\end{abstract}
\vspace{5mm}

\section{Introduction}

One way of solving explicitly a quantum integrable system is to
realize the system as a quantum reduction of a larger quantum
integrable system for which explicit solution may be easily found.
Choice of the larger system is obviously not unique and is a matter
of convenience. To practically implement this approach one should
represent wave functions of quantum integrable system as particular
matrix elements of suitable Lie algebra representations allowing
central characters. In this representation the reduction is realized
via proper  choice of the matrix element. Interesting examples arise
by considering  reductions  with respect to non-abelian Lie group
symmetries. Typical case (along with various kinds of the Calogero
models \cite{OP}) is given by families of quantum open  Toda chains
solved via special  matrix elements of principle series
representations of  non-abelian reductive Lie groups identified
 with Whittaker functions \cite{Ja},\cite{SCH},\cite{Ha} and \cite{Ko1},\cite{Ko2}.

Recall that another more traditional approach to integration of
quantum integrable system employs large abelian Lie group of
symmetries generated by mutually commuting integrable flows. The
large abelian symmetry allows  to find a proper set of quantum
canonical variables and to realize wave functions of the quantum
integrable systems as matrix elements of representations of the
Heisenberg Lie algebras. However, this straightforward  quantum
generalization of the classical integration algorithm encounters
various difficulties. The main reason is that quite non-trivial
realization of the corresponding abelian symmetries does not allow
to  fix properly operator ordering ambiguities in a simple and
explicit way.

It is natural to  try to combine these two complimentary approaches
by considering a larger explicitly integrable quantum theory
completely defined in terms of  its Lie symmetry given by extensions
of the Heisenberg Lie algebras. As it was demonstrated in the
previous short announcement \cite{GLO5} the theory of generalized
hypergeometric functions developed in \cite{GGZ}, \cite{GKZ}
supplies us with a large class of quantum integrable systems which
may be solved this way.

In this note we show that general GKZ hypergeometric function may be
identified with a matrix element in an irreducible representation of
a suitable multidimensional oscillator Lie algebra. The explicit
identification is most natural when is done in terms of
Gelfand-Graev hypergeometric functions
 (referred to as GG-functions). The   latter were introduced in
 \cite{GG} as  a more symmetric formulation   of the original
 GKZ hypergeometric functions.

The class of integrable systems solved in terms of GKZ hypergeometric functions
 includes open $\mathfrak{gl}_{\ell+1}$-Toda chains so that the corresponding
solutions given by $\mathfrak{gl}_{\ell+1}$-Whittaker functions
allow integral representations via GKZ integrals (see
\cite{BCFKvS}). The matrix element representation of  open
$\gl_{\ell+1}$-Toda chain wave functions introduced in this note
should be contrasted  with another one  based on representation
theory of reductive Lie algebra $\gl_{\ell+1}(\IR)$
\cite{Ko1},\cite{Ko2}. To stress the difference we provide a
detailed description of both formulations in the case of Toda type
models associated with maximal and minimal parabolic subalgebras of
$\gl_{\ell+1}$ \cite{BCFKvS},\cite{GLO4}. On a more fundamental
level the relation between two matrix element formulations will be
considered elsewhere.

Our interest in various realizations of the Whittaker functions
steams for the fact that the Whittaker functions play important role
in the formulation of local  Archimedean Langlands correspondence
(see \cite{GLO1},\cite{GLO3},\cite{G} and references therein). Thus
the identification of the Whittaker functions with matrix elements
of the non-reductive oscillator Lie algebra possibly leads to
interesting implications for Archimedean Langlands correspondence
related to various approaches based on abelianization/torification
and  mirror symmetry. This theme will be discussed in future
publications.

The plan of the paper is as follows. In Section 2 we provide basics
of GKZ hypergeometric structure and associated hypergeometric
functions. Following mostly \cite{GG} we also define GG-system as a
symmetric form of the GKZ  hypergeometric structure. The
corresponding GG-functions are  the main objects of interest in the
following Sections. In Section 3 matrix element realization of the
general GG hypergeometric functions via representations theory of
multidimensional oscillator algebras is presented. In Section 4 the
set of defining equations for GG hypergeometric functions is
rederived using matrix element representation introduced in Section
3. Finally in Section 5 the matrix element representation in terms
of oscillator Lie algebra and reductive Lie algebra
$\mathfrak{gl}_{\ell+1}$, for $\mathfrak{gl}_{\ell+1}$-Whittaker
functions associated with maximal and minimal parabolic subalgebras
are considered in detail.

{\it  Acknowledgements:} The research of D.R.L. was supported by RSF
grant 16-11-10075. The work of S.V.O. was partially supported by the
EPSRC grant EP/L000865/1.

\section{The GKZ and GG hypergeometric structures}

The  GKZ hypergeometric functions are defined as solutions to
system of differential equations associated with the
$(\CA,c)$-hypergeometric structures, first introduced and studied in
\cite{GGZ}, \cite{GKZ}. The GKZ hypergeometric structure is defined
by the following data. We chose an $N$-element subset
 \be\label{matA}
 \CA\,=\,\{a_1,a_2,\ldots, a_N\}\subset \IZ^m,
 \ee
such that $\CA$ generates $\IZ^m$ as abelian group. In addition, we
fix a complex vector,
 \be\label{C}
  c\,=\,(c^1,\ldots,c^m)\in \IC^m\,.
  \ee
The collection $\CA$ defines the $m\times N$-matrix of maximal
possible rank:
 \be\label{AAA}
  A\,=\,\|a^s_i\|\,\in\,\Mat_{m\times N}(\IZ)\,,\qquad\rk(A)=m
  \ee
where $ a_i=(a_i^1,\cdots , a_i^m)\in \IZ^m$ are elements of $\CA$.

Let $\mathbb{L}\subset \CA$ be the relations lattice of $\CA$:
 \be\label{L}
  \mathbb{L}\,=\,\{(l_1,\ldots, l_N)\in \IZ^N\,:\,l_1a_1+\cdots
  +l_Na_N=0\}\,.
 \ee
Choosing a basis $\{l^{\a},\,\a\in J\}$ of the lattice $\mathbb{L}$
indexed by $J=\{1,\ldots,N-m\}$:
 \be\label{LLL}
  \IL\,=\,\Span\bigl\{l^{\a}=(l_1^{\a},\,\ldots,\,l^{\a}_N)\,:\quad\a\in
  J\bigr\}\,,
  \ee
we consider the corresponding relation matrix
 \be\label{MMM}
  M\,=\,\|l^{\a}_i\|\,\in\,\Mat_{(N-m)\times N}(\IZ)\,,\qquad\rk(M)=N-m\,.
  \ee
The matrices $A$ and $M$ enjoy the orthogonality property
  \be\label{ORTH}
  AM^{\top}\,=\,0\,\in\,\Mat_{m\times (N-m)}(\IZ)\,,
 \ee
where $\top$ denotes the standard matrix transposition. The GKZ
hypergeometric function $f(u)$ is a solution to the GKZ-system
introduced below.
\begin{de}\label{GKZ1} The GKZ-system of differential equations
  associated with the data $(\CA,\,c)$ consists of the
  following set of  equations in variables
$u=(u_1,\cdots ,u_N)\in(\IR_+)^N$:
\begin{itemize}
\item[(1)]  For every $l\in \mathbb{L}$ one has
 \be\label{DEF}
  \prod_{i\in I\atop l_i<0} \left(-\frac{\pr}{\pr u_i}\right)^{-l_i}\,
  f\,=\,\prod_{i\in I\atop l_i>0} \left(-\frac{\pr}{\pr u_i}\right)^{l_i}\,
  f\,.
 \ee

\item[(2)] The  differential equations enumerated by $s\in\{1,\ldots,m\}$:
 \be\label{TORUS}
  a_1^su_1\frac{\pr f}{\pr u_1}+ \cdots+ a_N^su_N\frac{\pr f}{\pr
  u_N}=c^sf.
 \ee
\end{itemize}
\end{de}
The equations \eqref{TORUS} may be easily integrated and thus allow
a reduction of the solutions to the system of equations \eqref{DEF},
\eqref{TORUS} to functions in
  $(N-m)$-variables indexed by $J=\{1,\ldots,N-m\}$. However the
resulting equations on the function $(N-m)$-variables have more
complicated form.

Now we fix  a special solution to the GKZ-system \eqref{DEF},
\eqref{TORUS}.
\begin{prop}\label{GKZsol} Given the GKZ data $(\CA,c)$,  the GKZ-system
\eqref{DEF}, \eqref{TORUS} allows the following solution:
 \be\label{INT}
  f_{\g}(u)\,
  =\!\int\limits_{\IR_+^N}\!\prod_{i=1}^N\frac{dt_i}{t_i}\,
  t_i^{\g_i}\,e^{-u_it_i}\prod_{\a\in J}
  \delta\Big(\prod_{j\in I}t_j^{l^{\a}_j}\,-\,1\Big)\,,\qquad{\rm
  Re}(\g_i)>0\,.
 \ee
Here $\{l^\a,\,\a\in J\}\subset\IL$ is a basis \eqref{LLL} and
$\g=(\g_1,\ldots,\g_N)\in\IC^N$ is a vector subjected to
 \be\label{GKZsp}
  c^s+\sum_{j\in I}a^s_j\g_j\,=\,0\,,\qquad s\in\{1,\ldots,m\}\,.
  \ee
The function \eqref{INT} is independent of the choice of the basis
as well as of the choice of $\g\in\IC^N$ satisfying \eqref{GKZsp}.
\end{prop}

\proof Let $\g=(\g_1,\ldots,\g_N)\in\IC^N$ such that ${\rm
Re}(\g_i)>0,\,i\in I$. For every $l\in\IL$, let us verify the first
assertion by substituting \eqref{INT} into \eqref{DEF} and
\eqref{TORUS}. For the first equation,
 \be
  \prod_{i\in I\atop l_i>0}
  \Big(-\frac{\pr}{\pr u_i}\Big)^{l_i}
  \cdot f_{\g}(u)\,-  \prod_{i\in I\atop l_i<0}
  \Big(-\frac{\pr}{\pr u_i}\Big)^{-l_i}
  \cdot f_{\g}(u)\\
  =\!\int\limits_{\IR_+^N}\!\prod_{i=1}^N\frac{dt_i}{t_i}\,
  t_i^{\g_i}\,e^{-u_it_i}\,
  \Big(\prod_{i\in I\atop l_i>0}t_i^{l_i}-
  \prod_{i\in I\atop l_i<0}t_i^{-l_i}\Big)
    \prod_{\alpha\in J}  \delta\Big(\prod_{j\in I}t_j^{l^{\a}_j}\,-\,1\Big)
 \ee
 \be
   =\!\int\limits_{\IR_+^N}\!\prod_{i=1}^N\frac{dt_i}{t_i}\,
   t_i^{\g_i}\,e^{-u_it_i}\,
  \prod_{i\in I\atop l_i<0}t_i^{-l_i}\,\,\,
  \Big(\prod_{i\in I}t_i^{l_i}-1\Big)
    \prod_{\alpha\in J}  \delta\Big(\prod_{j\in
      I}t_j^{l^{\a}_j}\,-\,1\Big)=0.
 \ee
The last equality follows since for $l=\sum\limits_{\a\in
J}n_{\a}l^{\a}\in\IL$, we have
 \be
  \prod_{i\in I}t_i^{l_i}\,
  =\,\prod_{i\in I}t_i^{\sum\limits_{\a\in J}n_{\a}l^{\a}_i}\,
  =\,\prod_{\a\in J}\Big(\prod_{i\in
  I}t_i^{l^{\a}_i}\Big)^{n_{\a}}\,,
 \ee
which equals 1 by taking into account the delta-factors in the
integrand of \eqref{INT}.

To prove the second equation it is useful to  change the integration
variables $t_i\to u_i^{-1}t_i$, so the integral \eqref{INT} takes
the form:
 \be\label{INTn}
  f_{\g}(u)\,
  =\,\prod_{i\in I}u_i^{-\g_i}\!\int\limits_{\IR_+^N}\!
  \prod_{i=1}^N\frac{dt_i}{t_i}\,
  t_i^{\g_i}\,e^{-t_i}\prod_{\a\in J}
  \delta\Big(\prod_{j\in
  I}u_j^{-l^{\a}_j}t_j^{l^{\a}_j}\,-\,1\Big)\,.
 \ee
Then the equations \eqref{TORUS}
follow from the orthogonality relation \eqref{ORTH} between
$A$ and $M$.

For the solution \eqref{INT}, independence of the choice of the
basis $\{l^\a,\,\a\in J\}\subset\IL$ can be verified as follows. Let
$\{\lt^{\a},\,\a\in J\}\subset\IL$ be another basis and let
$\|g^{\a}_{\beta}\|\in GL_{N-m}(\IZ)$ be the transition matrix
between the two bases:
 \be
  \lt^{\a}\,=\,\sum_{\beta\in J}g^{\a}_{\beta}l^{\beta}\,, \qquad
  \det\|g^{\a}_{\beta}\|\in \IZ^*=\{\pm 1\}.
 \ee
By change of integration variables $T_i=\ln t_i,\,i\in I$,
\eqref{INT} takes the following form:
 \be\label{INTlog}
  f_{\g}(u)\,
  =\!\int\limits_{\IR^N}\!\prod_{i=1}^NdT_i\,
  e^{\g_iT_i\,-\,u_ie^{T_i}}\prod_{\a\in J}
  \delta\Big(\sum_{j\in I}l^{\a}_jT_j\Big)\,.
 \ee
For $T=(T_1,\ldots,T_N)\in\IR^N$, consider vector $S\in\IR^{N-m}$
with the coordinates $S^{\a}=\sum\limits_{j\in I}l^{\a}_jT_j$, and
define
 \be
 \delta(S)=\prod_{\a\in J} \delta(S^\a)=\prod_{\a\in J}
  \delta\Big(\sum_{j\in I}l^{\a}_jT_j\Big).
  \ee
Then applying standard rule for linear change of arguments in
delta-functions
 \be
\delta(gS)=\frac{1}{|\det(g)|}\,
 \delta(S)=\delta(S),\qquad\det(g)\in\IZ^*=\{\pm1\}\,,
 \ee
which verifies the independence of the choice of the basis
$\{l^\a\}\subset\IL$ for \eqref{INT}.

Finally, for $f_{\g}(u)$, independence of the choice of solution
$\gamma=(\g_1,\ldots,\g_N)\in\IC^N$ to the equations \eqref{GKZsp}
may be checked as follows. Given a pair of solutions,
 $\gamma$ and $\gamma'$, they are related by
 \be\label{SpecSym}
  \g_j=\g_j'+\sum_{\a\in J}l^\a_j\xi_\a, \qquad\xi=\|\xi_{\a}\|\in\IC^{N-m}\,.
 \ee
Then substituting this into the integral \eqref{INT} results in
$f_{\g}=f_{\g'}$ by taking into account the delta-factors in the
integrand. $\Box$

\begin{ex}\label{EX} Consider the two elementary examples of the
  GKZ-hypergeometric functions
associated with $\CA$, corresponding
    to the cases $|I|=1,\,|J|=0$ and $|I|=|J|=1$.
\begin{enumerate}
\item The  GKZ-hypergeometric function
  corresponding to the GKZ data with
$N=m=1$ and $J=\emptyset$ according to Definition \ref{GKZ1} is
given by
 \be\label{Gamma}
  f_{\g}(u)\,=\int\limits_{\IR_+}\!\frac{dt}{t}t^\gamma\,e^{-ut}\,
  =\,u^{-\gamma}\,\Gamma(\gamma)\,,\qquad{\rm
  Re}(\g)>0\,,
  \ee
and satisfies the following equation (a special case of \eqref{TORUS}):
 \be\label{GammaEq1}
  u\frac{\pr}{\pr u}f_{\g}(u)\,=\,-\g\,f_{\g}(u)\,.
 \ee
  Note that
the equation \eqref{DEF} is absent in this case.

\item The GKZ-hypergeometric function associated with the
data $|I|=|J|=1,\,m=0$ is given by the exponential function:
 \be\label{F}
  f_\g(u)\,
  =\,\int\limits_{\IR_+}\!\frac{dt}{t}\,t^{\g}\,e^{-ut}\,\delta(t-1)\,
  =\,e^{-u}\,,
 \ee
and satisfies the following equation:
 \be\label{Feq1}
  \Big(-\frac{\pr}{\pr u}\Big)f_\g(u)\,=\,f_\g(u)\,,
 \ee
which is a special case of \eqref{DEF}.  The equation
\eqref{TORUS} is absent in this case.
\end{enumerate}
\end{ex}

Definition \ref{GKZ1} of the GKZ system has an obvious asymmetry
that we would like to get rid of. Namely, it is natural to consider
a solution to \eqref{DEF},\eqref{TORUS} as a function both in
$u\in\IR_+^N$ and $c\in\IC^m$. Here the vector $c$ plays the role of
spectral variables. The difference between number of
$u$-variables and $c$-variables is taken into account by additional
symmetries generated by linear equations \eqref{TORUS}. Clearly
$u$-variables and $c$-variables are treated quite differently in
Definition \ref{GKZ1}.  On the other hand, the explicit solution
$f_{\g}(u)$ written down using spectral $\gamma$-variables of the
same number as $u$-variables. As a compensation, we gain an
additional symmetry over spectral parameters
$\g=(\g_1,\ldots,\g_N)\in\IC^N$
 reducing it effectively
 to the $m$ parameters $\{c^s,\,1\leq s\leq m\}$ via \eqref{SpecSym}.
 All this suggests a
completely symmetric definition of the GKZ system of differential
equations. Basic ingredients of this reformulation already appeared
in \cite{GG} and thus we will call this GG-hypergeometric structure.

Let the set $\CA$ be the same as in the definition of GKZ-data
$(\CA,c)$ and let $\IL$ be the corresponding relation lattice
\eqref{L}. Let us introduce GG-hypergeometric system associated with
$\CA$. The GG-hypergeometric function $\Phi_{\g}(u)$ is a solution
to the GG-system.

\begin{de}\label{DD2} The  GG-system associated with $\CA$ consists of the
  following set of  equations in variables
  $u=(u_1,\ldots ,u_N)\in\IR_+^N$  and $\g=(\g_1,\ldots ,\g_N)\in \IC^N$:
\begin{itemize}
\item[(1)]  For every $l\in \mathbb{L}$ one has
 \be\label{DEF1}
 \prod_{i\in I\atop l_i<0}
 \left(-\frac{\pr}{\pr u_i}+\frac{\g_i}{u_i}\right)^{-l_i}\,
 \Phi\,=\,\prod_{i\in I\atop l_i>0}
 \left(-\frac{\pr}{\pr u_i}+\frac{\g_i}{u_i}\right)^{l_i}\,
  \Phi\,.
 \ee

\item[(2)]  For every $l\in \mathbb{L}$ one has
\be\label{DEF1d2}
  \sum_{i\in I}l_i \left(\frac{\pr}{\pr\g_i}-\ln u_i\right)\,\,\Phi\,=0\,,
 \ee
or, equivalently,
   \be\label{DEF1d}
  \prod_{i\in I}e^{\theta\lambda_i\pr_{\g_j}}\Phi\,
  =\,\prod_{i\in I}u_i^{\theta\lambda_i}\,\Phi\,,\qquad \forall\theta \in
  \IR\,.
 \ee

\item[(3)] The system of $m$ differential equations for $s\in\{1,\ldots,m\}$:
 \be\label{TORUS1}
  \Big\{a_1^su_1\frac{\pr }{\pr u_1}+ \cdots+ a_N^su_N\frac{\pr}{\pr
  u_N}\Big\}\Phi\,=\,0\,.
 \ee

\item[(4)] The system of $m$ difference equations for $s\in\{1,\ldots,m\}$:
 \be\label{TORUS1d}
  \Big\{a^s_1\left(e^{\frac{\pr}{\pr\g_1}}-\gamma_1\right)+\ldots
  +a_N^s\left(e^{\frac{\pr}{\pr\g_N}}-\gamma_N\right)\Big\}\Phi\,=\,0\,.
 \ee
\end{itemize}
\end{de}

Note that in comparison with Definition \ref{GKZ1} we introduce
additional dual difference equations over the spectral variables
$\gamma\in\IC^N$.

To contrast Definitions \ref{GKZ1} and \ref{DD2} it is useful to
consider the two simple cases  continuing Example \ref{EX} above.

\begin{ex}
\begin{enumerate}
\item The  GG-hypergeometric function
  corresponding to $\CA$  with
$N=m=1$ and $J=\emptyset$ according to Definition \ref{DD2} is given
by the  standard Gamma-function:
 \be
  \Phi_\gamma(u)=\int\limits_{\IR_+}\!\frac{dt}{t}\,t^{\g}\,e^{-t}=\Gamma(\gamma)\,,\qquad{\rm
  Re}(\g)>0\,,
 \ee
satisfying the equations
 \be\label{GammaEq1d}
  \frac{\pr}{\pr u}\,\Phi_\gamma(u)\,=\,0\,,
 \ee
and
 \be\label{GammaEq2d}
  \left(e^{\pr_\gamma}-\gamma\right)\,\Phi_\gamma(u)\,=\,0\,.
 \ee
Here \eqref{GammaEq1d} is an instance of \eqref{TORUS1}, meanwhile
\eqref{GammaEq2d} is an instance of \eqref{TORUS1d}. Note that the
equations \eqref{DEF1} and \eqref{DEF1d2},\eqref{DEF1d} are absent
in this case.

\item The GG-hypergeometric function associated with $\CA$ such that $|I|=|J|=1,\,m=0$ is given by
 \be
 \Phi_{\g}(u)\,=\,u^{\g}\,e^{-u}\,,
 \ee
satisfying the differential equation:
 \be\label{Feq1d}
  \Big\{u\frac{\pr}{\pr u}\,-\,\g\Big\}\Phi_\gamma(u)\,=\,u\,\Phi_\gamma(u)\,,
 \ee
and the difference equation:
 \be\label{Feq2d}
  \Big\{e^{\pr_{\g}}-u\Big\}\Phi_\gamma(u)\,=\,0\,.
 \ee
Here \eqref{Feq1d} is an instance of \eqref{DEF1}, meanwhile
\eqref{Feq2d} is an instance of \eqref{DEF1d}. The equations
\eqref{TORUS1},\eqref{TORUS1d} are absent in this case.
\end{enumerate}
\end{ex}

A relation between the two formulations of the GKZ-hypergeometric
structures is manifested by the following.

\begin{prop} Given the GKZ-datum  $\CA$,
the GG-system \eqref{DEF1}-\eqref{TORUS1d}   allows the following
solution:
 \be\label{INT1}
  \Phi_{\g}(u)\,
  =\!\int\limits_{\IR_+^N}\!\prod_{i=1}^N\frac{dt_i}{t_i}\,
  t_i^{\g_i}\,e^{-t_i}\prod_{\a\in J}
  \delta\Big(\prod_{j\in
  I}u_j^{-l^{\a}_j}t_j^{l^{\a}_j}\,-\,1\Big)\,,\qquad{\rm
  Re}(\g_i)>0\,,
 \ee
where  $\{l^\a,\,\a\in J\}\subset\IL$ is the basis  \eqref{LLL}. The
solution \eqref{INT1} is independent of the choice of basis
$\{l^\a,\,\a\in J\}\subset\IL$ and may be identified with
\eqref{INT} via
 \be\label{INTid}
  \Phi_{\g}(u)\,=\,\Big(\prod_{i\in I}u_i^{\g_i}\Big) \times f_{\g}(u)\,.
 \ee
\end{prop}
\proof Let $\g\in\IC^N$ be such that ${\rm
  Re}(\g_i)>0,\,i\in I$. Similarly to the proof of
Proposition \ref{GKZsol} above, one checks that \eqref{INT1}
satisfies \eqref{DEF1}, \eqref{DEF1d2} and \eqref{TORUS1} in
straightforward way. For the integral \eqref{INT1}, independence of
the basis choice follows from \eqref{INTid} and from the
independence for $f_{\g}(u)$ by Proposition \ref{GKZsol}.

One observes that the equations \eqref{DEF1} and \eqref{TORUS1} are
identified with \eqref{DEF} and \eqref{TORUS} via \eqref{INTid} and
the following identity for every $i\in I$:
 \be
  \prod_{j\in I}u_j^{\g_j}\Big(\frac{\pr}{\pr u_i}\Big)^n\prod_{j\in I}u_j^{-\g_j}\,
  =\,\Big(\frac{\pr}{\pr u_i}\,-\,\frac{\g_i}{u_i}\Big)^n\,,\qquad n\in\IZ_{\geq0}\,.
 \ee
To verify \eqref{TORUS1d}, consider the integrand of \eqref{INT1}:
 \be
  F_{\g}(t)\,
  =\,\prod_{i\in I}t_i^{\g_i}e^{-t_i}
  \prod_{\a\in J}\delta\Big(\prod_{j\in I}u_j^{-l^{\a}_j}t_j^{l^{\a}_j}\,-\,1\Big)\,.
 \ee
Then its differential is given by
 \be
  dF_{\g}(t)\,=\,\sum_{i\in I}\frac{\pr F_{\g}(t)}{\pr t_i}\,dt_i\,
  =\,\sum_{i\in I}\bigl(\g_i-t_i\bigr)F_{\g}(t)\frac{dt_i}{t_i}\\
  +\,\sum_{i\in I}\sum_{\a\in J}l^{\a}_i
  \prod_{j\in I}u_j^{-l^{\a}_j}t_j^{l^{\a}_j}
  \delta'\Big(\prod_{j\in I}u_j^{-l^{\a}_j}t_j^{l^{\a}_j}\,-\,1\Big)
  \prod_{\beta\in J\atop\beta\neq\a}
  \delta\Big(\prod_{j\in
  I}u_j^{-l^{\beta}_j}t_j^{l^{\beta}_j}\,-\,1\Big)\frac{dt_i}{t_i}\,.
 \ee
Therefore, observing that
 \be
  \bigl(\g_i-t_i\bigr)F_{\g}(t)\,=\,\Big(\g_i-e^{\pr_{\g_i}}\Big)F_{\g}(t)\,,
 \ee
by the orthogonality relation \eqref{ORTH}, one deduces the
following, for $1\leq s\leq k+1$:
 \be
  \sum_{i\in I}a^s_i\Big(\g_i-e^{\pr_{\g_i}}\Big)\Phi_{\g}(u)\,
  =\!\int\limits_{\IR_+^N}\!\prod_{j=1}^N\frac{dt_j}{t_j}
  \sum_{i\in I}a^s_i\bigl(\g_i-t_i\bigr)F_{\g}(t)\\
  =\sum_{i\in I}a^s_i\!\int\limits_{\IR_+^{N-1}}\!\prod_{j\in I\atop j\neq i}
  \frac{dt_j}{t_j}\int\limits_{\IR_+}dt_i\,\frac{\pr F_{\g}(t)}{\pr
  t_i}\,=\,0\,,
 \ee
since $F_{\g}(t)\big|_{t_i=0}=F_{\g}(t)\big|_{t_i=+\infty}=0$ for
each $i\in I$.\, $\Box$

As it is already mentioned in Introduction, the system of equations
\eqref{DEF1}-\eqref{TORUS1d} is naturally related to the GG-system
introduced and studied in \cite{GG}. In particular, the integral
\eqref{INT1} can be identified with the GG-function as it is defined
in \cite{GG}. Alternatively, in \cite{GG} the function \eqref{INT1}
is described as a solution to \eqref{DEF1}, \eqref{DEF1d2},
\eqref{TORUS1} and the additional set of first-order
differential-difference equations \eqref{DDeq} from the following
Proposition.

\begin{prop} The GG-function \eqref{INT1} satisfies
  the following system of  GG-equations:
 \be\label{DDeq}
  \Big\{-u_i\frac{\pr}{\pr u_i}+\g_i\Big\}\cdot\Phi_{\g}(u)\,
  =\,e^{\pr_{\g_i}}\cdot\Phi_{\g}(u),\qquad i\in I\,.
 \ee
\end{prop}

\proof To verify the assertion, for every $j\in I$, we substitute
\eqref{INT1} (after changing $t_i\to u_it_i,\,i\in I$) into
\eqref{DDeq}:
 \be
  \Big\{-u_j\frac{\pr}{\pr u_j}+\g_j\Big\}\Phi_{\g}(u)\,
  =\!\int\limits_{\IR_+^N}\!\prod_{i=1}^N\frac{dt_i}{t_i}\,\prod_{\a\in J}
  \delta\Big(\prod_{j\in I}t_j^{l^{\a}_j}\,-\,1\Big)\\
  \times\Big\{-u_j\frac{\pr}{\pr u_j}+\g_j\Big\}
  \prod_{i=1}^N(u_it_i)^{\g_i}\,e^{-u_it_i}\\
  =\,\!\int\limits_{\IR_+^N}\!\prod_{i=1}^N\frac{dt_i}{t_i}\,\prod_{\a\in J}
  \delta\Big(\prod_{j\in I}t_j^{l^{\a}_j}\,-\,1\Big)\,
  (u_jt_j)\prod_{i=1}^N(u_it_i)^{\g_i}\,e^{-u_it_i}\\
  =\,\Phi_{\g_1,\ldots,\g_j+1,\ldots,\g_N}(u)\,
  =\,e^{\pr_{\g_j}}\cdot\Phi_{\g}(u)\,.
 \ee
 $\Box$

One might note that \eqref{DDeq} and \eqref{TORUS1} entail the dual
difference equations \eqref{TORUS1d}.

\section{Matrix element representation}

In this section we provide a representation of the GG-hypergeometric
functions \eqref{INT1} as matrix elements in  irreducible
representation of a suitable class of non-reductive Lie algebras.
Namely, for $I=\{1,\ldots,N\}$, let $\fL_N$ be the Lie algebra
generated by the central element $\CC$ and by
$\{\CE_i,\CF_i,\CH_i,\,i\in I\}$ subjected to the following defining
relations:
 \be\label{LA}
  [\CE_i,\CF_j]=\CC\delta_{ij}\,, \qquad[\CH_i,\CE_j]=-\delta_{ij}\CE_j\,,\qquad
  [\CH_i,\CF_j]=\delta_{ij}\CF_j\,, \qquad i,j\in I\,.
 \ee
The subalgebra in $\fL_N$ generated by $\CC$ and by
$\{\CE_i,\,\CF_i,\,i\in I\}$ is isomorphic to the
$(2N+1)$-dimensional Heisenberg algebra. Thus the Lie algebra
$\fL_N$ may be considered as multidimensional version of the
standard oscillator Lie algebra.

For $\g=(\g_1,\ldots,\g_N)\in\IC^N$ such that ${\rm
Re}(\g_i)>0,\,i\in I$, let $\pi_{\g}$ be the $\fL_N$-representation
in the Schwartz space $\CV_{\g}$ of smooth functions in
$t=(t_1,\ldots,t_N)\in(\IR_+)^N$ decreasing rapidly with all its
derivatives at infinity:
 \be
  \CV_{\g}=\bigl\{f\in C^{\infty}(\IR_+^N):
  \lim_{t_i\to\infty}\!\bigl(t_1^{n_1}\cdots t_N^{n_N}|\pr_{t_1}^{m_1}\cdots\pr_{t_N}^{m_N}\!f(t)|\bigr)=0,\,
  \forall n_i,m_i\in\IZ_{\geq0}\bigr\}.
 \ee
Namely, the representation $(\pi_{\g},\CV_{\g})$ is defined by the
following action of the generators \eqref{LA}:
 \be\label{rep}
  \pi_{\g}(\CC)=1,\quad \pi_{\g}(\CE_i)=-\pr_{t_i}\,,\quad
  \pi_{\g}(\CF_i)=-t_i\,,\quad
  \pi_{\g}(\CH_i)=\g_i+t_i\pr_{t_i}\,,
  \quad i\in I\,.
 \ee
Let $\CU(\fL_N)$ be the universal enveloping algebra of $\fL_N$,
 then $\CV_{\g}$ affords a structure of
$\CU(\fL_N)$-module. Moreover, the action of generator
$\CH_i\in\fL_N$ in $\CV_{\g}$ integrates to the action of the group
$(\IR_+)^N$. It will be  also important in the following that the
generators $\CF_i,\,i\in I$  are invertible in the representation
$\CV_{\g}$, therefore $(\pi_{\g},\CV_{\g})$ extends to
representation of a larger algebra containing negative powers of
$\CF_i,\,i\in I$.  We fix a dual module $\CV^{\vee}_{\g}$ realized
in the space of generalized functions on $\IR_+^N$ and denote the
corresponding contragradient representation by $\pi^{\vee}_{\g}$.
There is a non-degenerate $\fL_N$-invariant pairing
$\<\,\,,\,\>\,:\CV^{\vee}_{\g}\times\CV_{\g}\to\IC$,
 \be\label{pairing}
  \<\phi,\,\varphi\>=\int\limits_{\IR_+^N}\prod_{i\in I}dt\,
  \phi(t)\,\varphi(t)\,,\qquad\varphi\in\CV_{\g}\,,\quad\phi\in\CV^{\vee}_{\g}\,,\\
  \<\pi^{\vee}_{\g}(X)\cdot\phi,\,\varphi\>\,
  =\,-\<\phi,\,\pi_{\g}(X)\cdot\varphi\>,\qquad\forall
  X\in\fL_N\,.
 \ee
Then the $\fL_N$-action in the dual module $\CV^{\vee}_{\g}$ can be
computed explicitly:
 \be\label{dualrep}
  \pi^{\vee}_{\g}(\CC)=-1,\qquad\pi^{\vee}_{\g}(\CE_i)=-\pr_{t_i}\,,\qquad
  \pi^{\vee}_{\g}(\CF_i)=t_i\,,\\
  \pi^{\vee}_{\g}(\CH_i)=1-\g_i+t_i\pr_{t_i}\,,\qquad i\in I\,.
 \ee

Given the GKZ-datum $\CA$, let us  fix
 a basis $\{l^{\a},\,\a\in J\}$ in the relation lattice
 $\IL$ \eqref{LLL}. Let  us introduce a vector
$\phi_R\in\CV_{\g}$ and a covector $\phi_L\in\CV^{\vee}_{\g}$
defined by
 \be\label{psiR}
  \pi_{\g}(\CE_i)\cdot\phi_R\,=\,\phi_R,\qquad i\in I\,,
 \ee
and
 \be\label{psiL}
  \prod_{i\in I\atop l^{\a}_i<0}\pi^{\vee}_{\g}(\CF_i)^{-l^{\a}_i}\cdot\phi_L\,
  =\,\prod_{i\in I\atop l^{\a}_i>0}
  \pi^{\vee}_{\g}(\CF_i)^{l^{\a}_i}\cdot\phi_L\,,\quad
  \a\in\{1,\ldots,N-m\}\,,\\
  \sum_{i=1}^Na^s_i\,\pi^{\vee}_{\g}(\CH_i)\cdot\phi_L\,=\,0\,,\qquad
  s\in\{1,\ldots,m\}\,.
  \ee
Note that invertibility of the images of the generators $\CF_i$ in
the considered representation allows to rewrite equivalently the
first line of \eqref{psiL} as follows:
 \be
  \prod_{i\in I}\pi^{\vee}_{\g}(\CF_i)^{l^{\a}_i}\cdot\phi_L\,
  =\,\phi_L\,,\qquad\a\in\{1,\ldots,N-m\}\,.
 \ee
Furthermore, for an arbitrary $l=\sum\limits_{\a\in
J}n_{\a}l^{\a}=(l_1,\cdots ,l_N) \in\IL$, the following holds:
 \be\label{psiLext}
  \prod_{i\in I}\pi^{\vee}_{\g}(\CF_i)^{l_i}\cdot\phi_L\,
  =\,\prod_{i\in I}\prod_{\a\in J}\pi^{\vee}_{\g}(\CF_i)^{n_{\a}l^{\a}_i}\cdot\phi_L\\
  =\,\prod_{\a\in J}
  \Big(\prod_{i\in I}\pi^{\vee}_{\g}(\CF_i)^{l^{\a}_i}\Big)^{n_{\a}}\cdot\phi_L\,
  =\,\phi_L\,.
 \ee

\begin{lem} For $\phi_L,\,\phi_R$
defined in \eqref{psiL}, \eqref{psiR}, the following explicit
expressions hold
 \be\label{psiLsol}
  \phi_R(t)=e^{-\sum\limits_{i=1}^N t_i},\qquad
  \<\phi_L,\varphi\>
  =\!\int\limits_{\IR_+^N}\!\prod_{i\in I}\frac{dt_i}{t_i}\,\,
  \prod_{i\in I}t_i^{\g_i}
  \prod_{\a\in J}\delta\Big(\prod_{j\in
  I}t_j^{l^{\a}_j}\,-\,1\Big)\,\varphi(t)\,,
 \ee
for an arbitrary $\varphi\in\CV_{\g}$.
\end{lem}
\proof The solution to \eqref{psiR} is given by the first expression
in \eqref{psiLsol}. As for $\phi_L$, the first defining relation in
\eqref{psiL} is a direct   consequence of presence of the
delta-factors in expression for $\phi_L$ in \eqref{psiLsol}. For the
second relation, substituting $\pi^{\vee}_{\g}(\CH_j)$ from
\eqref{dualrep} into \eqref{psiLsol}, for each $s\in\{1,\ldots,m\}$
implies
 \be
  \sum_{j=1}^Na^s_j\<\pi^{\vee}_{\g}(\CH_j)\cdot\phi_L,\varphi\>\\
  =\sum_{j=1}^Na^s_i\!
  \int\limits_{\IR_+^N}\!\prod_{i\in I}dt_i\,\varphi(t)\,
  (1-\g_j+t_j\pr_{t_j})\prod_{i\in I}t_i^{\g_i-1}
  \prod_{\a\in J}\delta\Big(\prod_{j\in I}t_j^{l^{\a}_j}\,-\,1\Big)\\
  =\,\sum_{\a\in J}\sum_{j=1}^Na^s_jl^{\a}_j
  \int\limits_{\IR_+^N}\!\prod_{i\in I}\frac{dt_i}{t_i}\,\,
  t_i^{\g_i\,+\,l^{\a}_i}
  \delta'\Big(\prod_{j\in I}t_j^{l^{\a}_j}\,-\,1\Big)
  \prod_{\beta\neq\a}\delta\Big(\prod_{j\in
    I}t_j^{l^{\beta}_j}\,-\,1\Big)\,\varphi(t)\,,
 \ee
which vanishes due to the orthogonality relation
$\sum\limits_{i=1}^Na^s_il^{\a}_i=0$ in \eqref{ORTH}. $\Box$

Let us note that the conditions \eqref{psiL} on the covector
$\phi_L$  may be linearized as follows. Consider the subalgebra in
$\fL_N$ generated by $\{\CH_i,\,\CF_i,\,i\in I\}$, subjected to
  \be
  \CH_i\CF_j-\CF_j\CH_i=\delta_{ij} \CF_j\,, \qquad i,j\in I\,.
  \ee
Then it may be embedded into the appropriately completed Heisenberg
algebra $\fH_N$ generated by the $P_i,Q_i,\,i\in I$ subjected to
 \be
 P_iQ_j-Q_jP_i=\delta_{ij}\,, \qquad i,j\in I\,.
 \ee
Explicitly, we define the following embedding
 \be
  \CH_i\longmapsto P_i\,, \qquad \CF_i\longmapsto e^{Q_i}\,.
 \ee
Then the defining relations \eqref{psiL} for covector $\phi_L$ are
equivalent to the condition of annihilation of $\phi_L$ by the
following $|I|=N$ operators:
 \be
  A_s=\sum_{i=1}^Na_i^sP_i\,,\quad s\in\{1,\ldots,m\};\qquad
  B_{\a}=\sum_{i=1}^Nl^{\a}_iQ_i\,,\quad\a\in J\,.
 \ee
These operators generate an $N$-dimensional commutative subalgebra
in the Heisenberg algebra $\fH_N$ and define a linear polarization
that in general differs from the one defined by maximal commutative
subalgebra generated by $P_i,\,i\in I$.

\begin{te} The integral solution \eqref{INT1} to the GG-system \eqref{DEF1}-\eqref{TORUS1d}
allows the following expression in terms of the matrix element in
the $\CU(\fL_N)$-representation $(\pi_{\g},\CV_{\g})$ for
$\g=(\g_1,\ldots,\g_N)\in\IC^N,\,{\rm Re}(\g_i)>0$:
 \be\label{wavefun}
  \Phi_{\g}(e^{y_1},\ldots,e^{y_N})
  =\,
  \bigl\<\phi_L,\,\pi_{\g}\Big(e^{\sum\limits_{j=1}^Ny_j\CH_j}\Big)\phi_R\bigr\>\,,
 \ee
where $\phi_R$ and $\phi_L$ are defined in \eqref{psiR} and
\eqref{psiL}, respectively.
\end{te}
\proof Substituting \eqref{psiLsol} into the matrix element we have
 \be
  \Phi_{\g_1,\ldots,\g_N}(e^{y_1},\ldots,e^{y_N})
  =\bigl\<\phi_L\,,\pi_{\g}\Big(e^{\sum\limits_{j=1}^Ny_j\CH_j}\Big)\phi_R\bigr\>\\
  =\int\limits_{(\IR_+)^N}\!\prod_{i=1}^N\frac{dt_i}{t_i}\,
  \phi_L(t)\,e^{\sum\limits_{j=1}^Ny_j\{\g_i+t_j\pr_{t_j}\}}\phi_R(t)\\
  =\,\int\limits_{(\IR_+)^N}\!\prod_{i=1}^N\frac{dt_i}{t_i}\,
  t_i^{\g_i}\,e^{\g_iy_i-t_ie^{y_i}}
  \prod_{\a\in J}\delta\Big(\prod_{j\in I}
  t_j^{l^{\a}_j}\,-\,1\Big)\,,
 \ee
which coincides with \eqref{INT1} for $u_i=e^{y_i}$ via substitution
$t_i\mapsto u_i^{-1}t_i,\,i\in I$. $\Box$

\section{Differential equations satisfied by the GKZ functions}

In this Section we derive the defining equations \eqref{DEF1},
\eqref{TORUS1}  using matrix element representation \eqref{wavefun}.
The way to derive these equations is similar to the way how
equations on matrix elements  arise from the action of Casimir
elements. Note that the dual equations
  \eqref{DEF1d2}-\eqref{DEF1d} and \eqref{TORUS1d} include
  differentiation over spectral parameters and arise in a different way
  by considering intertwining operators acting between different
  representations of oscillator algebra $\fL_N$.

Now for $\g=(\g_1,\ldots,\g_N)\in\IC^N$, given the
$\CU(\fL_N)$-module $(\pi_{\g},\,\CV_{\g})$ \eqref{rep}, let
$\CI_{\g}\subset\CU(\fL_N)$ be the primitive annihilation ideal:
 \be
  \CI_{\g}\,
  =\,\bigl\{X\in\CU(\fL_N)\,:\,\pi_{\g}(X)=0\bigr\}\,.
 \ee
For every $l=(l_1,\ldots,l_N)\in\IZ^N$, introduce the following
element in $\CU(\fL_N)$:
 \be\label{Cas}
  \CC(l)
  =\prod_{i\in I\atop l_i<0}\CF_i^{-l_i}\CE_i^{-l_i}\!
  \prod_{i\in I\atop l_i>0}\prod_{k=0}^{l_i-1}(\CH_i-\g_i-k)\,
  -\,\prod_{i\in I\atop l_i>0}\CF_i^{l_i}\CE_i^{l_i}\!
  \prod_{i\in I\atop l_i<0}\prod_{k=0}^{|l_i|-1}\!(\CH_i-\g_i-k)\,.
 \ee
\begin{lem}\label{Casimir0} For an arbitrary $l\in\IZ^N$, the element $\CC(l)$
belongs to the ideal $\CI_{\g}$.
\end{lem}
\proof Substitution of \eqref{rep} into \eqref{C} reads
 \be\label{CasimirRep}
  \pi_{\g}\bigl(\CC(l)\bigr)\,
  =\,\prod_{i\in I\atop l_i<0}t_i^{-l_i}(\pr_{t_i})^{-l_i}
  \prod_{i\in I\atop l_i>0}\prod_{k=0}^{l_i-1}(t_i\pr_{t_i}-k)\\
  -\,\prod_{i\in I\atop l_i>0}t_i^{l_i}(\pr_{t_i})^{l_i}
  \prod_{i\in I\atop
  l_i<0}\prod_{k=0}^{-l_i-1}(t_i\pr_{t_i}-k)\,.
 \ee
Then by the following identity
 \be\label{ID1}
  \prod_{k=0}^{n-1}(t\pr_t-k)\,=\,t^n(\pr_t)^n\,,\qquad\forall n>0\,,
 \ee
each of the two terms in \eqref{CasimirRep} equals $\prod_{i\in
I}(t_i)^{|l_i|}(\pr_{t_i})^{|l_i|}$ and hence cancel each other.
$\Box$

Now given the GKZ-datum  $\CA$, let $\IL$ be the corresponding
relation lattice.
\begin{prop} For each $l\in\IL$, the fact that
  $\pi_{\g}(\CC(l))=0$ entails the GG-equations \eqref{DEF1} satisfied by
the matrix element \eqref{wavefun}:
 \be\label{GKZeq}
 \Big\{\prod_{i\in I\atop l_i<0}\Big(-\frac{\pr}{\pr u_i}\,+\,
 \frac{\g_i}{u_i}\Big)^{-l_i}\,
  -\,\prod_{i\in I\atop l_i>0}\Big(-\frac{\pr}{\pr u_i}\,+\,
  \frac{\g_i}{u_i}\Big)^{l_i} \Big\}\,\Phi_{\g}(u)\,=\,0\,.
 \ee
\end{prop}

\proof Using the invertibility of $\pi_{\g}(\CF_i),\,i\in I$
   in the $\CU(\fL_N)$-representation $(\pi_{\g},\CV_{\g})$,  we
 introduce the elements $\widetilde{\CC}(l)\in\CU(\fL_N)$ defined by
 \be\label{Casimir}
 \pi_{\g}(\widetilde{\CC}(l))\,=\,
 \prod_{i\in I\atop l^{\a}_i<0}\pi_\g(\CF_i)^{l^{\a}_i} \times
  \pi_{\g}(\CC(l)),\,\qquad l\in\IL\,,
 \ee
also acting  by zero in $\CV_\g$. Substituting the factors
$\prod\limits_{k=0}^{n-1}(\CH_i-\g_i-k)$ from \eqref{Cas} into the
matrix element \eqref{wavefun} gives
 \be
  \bigl\<\phi_L,\,
  \pi_{\g}\Big(\prod_{k=0}^{n-1}(\CH_i-\g_i-k)\,e^{\sum_jy_j\CH_j}\Big)\phi_R\bigr\>\,
  =\,\prod_{k=0}^{n-1}(\pr_{y_i}-\g_i-k)\Phi_{\g}(e^y)\,.
 \ee
Insertion of \eqref{Casimir} into matrix element
 \eqref{wavefun}  and taking into account defining relations on the
right and left vectors \eqref{psiR} and
\eqref{psiL},\eqref{psiLext}, implies the following:
 \be
  \bigl\<\phi_L,\,
  \pi_{\g}\Big(\widetilde{\CC}(l)\,e^{\sum_jy_j\CH_j}\Big)\phi_R\bigr\>\\
   =\,\prod_{i\in I\atop l_i>0}\prod_{k=0}^{l_i-1}\Big(u_i\frac{\pr}{\pr u_i}-\g_i-k\Big)\,
  \bigl\<\phi_L,\,
  \pi_{\g}\Big(\prod_{i\in I\atop l_i<0}\CE_i^{-l_i}\,
  e^{\sum_jy_j\CH_j}\Big)\phi_R\bigr\>\\
  -\,\prod_{i\in I\atop  l_i<0}\prod_{k=0}^{-l_i-1}
  \Big(u_i\frac{\pr}{\pr u_i}-\g_i-k\Big)\,
  \bigl\<\phi_L,\,
  \prod_{i\in I}\pi_{\g}(\CF_i)^{l_i}\, \pi_{\gamma}\Big(
  \prod_{i\in I\atop l_i>0}\CE_i^{l_i}\,
  e^{\sum_jy_j\CH_j}\Big)\phi_R\bigr\>
 \ee
 \be
  =\,\prod_{i\in I\atop
    l_i>0}\prod_{k=0}^{l_i-1}\Big(u_i\frac{\pr}{\pr u_i}-  \g_i-k\Big)\,
  \prod_{i\in I\atop l_i<0}u_i^{-l_i}\,
  \bigl\<\phi_L,\,\pi_{\g}\Big(e^{\sum_jy_j\CH_j}\Big)\phi_R\bigr\>\\
  -\prod_{i\in I\atop l_i<0}\!\!\prod_{k=0}^{-l_i-1}\!\!\Big(u_i\frac{\pr}{\pr
   u_i}-\g_i-k\Big)\!\!\prod_{i\in I\atop l_i>0}u_i^{l_i}
  \bigl\<\prod_{i\in I}(-1)^{l_i}\pi^{\vee}_{\g}(\CF_i)^{l_i}\phi_L,\pi_{\g}\Big(
   e^{\sum_jy_j\CH_j}\Big)\phi_R\bigr\>\,.
 \ee
Then applying the relations \eqref{psiLext}, this results in
 \be
  \bigl\<\phi_L,\,
  \pi_{\g}\Big(\widetilde{\CC}(l)\,e^{\sum_jy_j\CH_j}\Big)\phi_R\bigr\>
  =\,\Big\{\prod_{i\in I\atop l_i>0}\prod_{k=0}^{l_i-1}\Big(u_i\frac{\pr}{\pr u_i}-\g_i-k\Big)\,
  \prod_{i\in I\atop l_i<0}u_i^{-l_i}\\
  -\,\prod_{i\in I}(-1)^{l_i}\prod_{i\in I\atop l_i<0}\prod_{k=0}^{-l_i-1}\Big(u_i\frac{\pr}{\pr
   u_i}-\g_i-k\Big)\,\prod_{i\in I\atop l_i>0}u_i^{l_i}\Big\}\,\Phi_{\g}(u)\\
   =\,\prod_{i\in I\atop  l_i>0}(-1)^{l_i}\prod_{i\in I}u_i^{|l_i|}
   \Big\{
  \prod_{i\in I\atop  l_i>0}(-u_i)^{-l_i}\prod_{k=0}^{l_i-1}\Big(u_i\frac{\pr}{\pr u_i}-\g_i-k\Big)\\
  -\,\prod_{i\in I\atop l_i<0}(-u_i)^{l_i}
  \prod_{k=0}^{-l_i-1}\Big(u_i\frac{\pr}{\pr
  u_i}-\g_i-k\Big)\Big\}\,\Phi_{\g}(u)\\
   =\,\prod_{i\in I\atop  l_i>0}(-1)^{l_i}
   \prod_{i\in I}u_i^{|l_i|}\Big\{\prod_{i\in I\atop  l_i>0}
  \Big(-\frac{\pr}{\pr u_i}+\frac{\g_i}{u_i}\Big)^{l_i}\,
  -\,\prod_{i\in I\atop l_i<0}\Big(-\frac{\pr}{\pr u_i}+
  \frac{\g_i}{u_i}\Big)^{-l_i}\Big\}\,\Phi_{\g}(u)\,.
 \ee
where in the last line \eqref{ID1} is used. Thus equating the above
expression to zero and dividing by the invertible function
$\prod\limits_{i\in I\atop  l_i>0}(-1)^{l_i}
   \prod\limits_{i\in I}u_i^{|l_i|}$ yields \eqref{GKZeq}. $\Box$

\section{Whittaker functions via  GKZ structures}

Investigations of quantum cohomology of partial flag manifolds give
rise to new families of quantum integrable systems of Toda type. The
corresponding wave functions (generating functions) are given by
generalized Whittaker functions associated with choice of parabolic
subalgebras $\frak{p}\subset\gl_{\ell+1}$; for the motivated
discussions of these see \cite{Giv},\cite{BCFKvS} and
\cite{GLO2},\cite{GLO4},\cite{O}. As in the original construction of
the standard Whittaker function (see \cite{Ja},\cite{SCH},\cite{Ha}
and \cite{Ko1},\cite{Ko2}), the generalized Whittaker functions can
be presumably realized as particular matrix elements of the
principal  series $\CU(\gl_{\ell+1})$-representations.  One also
expects that the generalized Whittaker functions allow expressions
in terms the appropriate GKZ  hypergeometric functions. Below we
consider two special instances of the $\gl_{\ell+1}$-Whittaker
functions associated with minimal and maximal parabolic subalgebras
and provide their expressions in terms of GKZ hypergeometric
functions. Various parts of these results may be found in
  references mentioned above and are included here for completeness.

Let us first fix the following notations for the general linear Lie
algebra and its subalgebras.  Let $V=\IC^{\ell+1}$ be a $\IC$-vector
space and let $\gl_{\ell+1}=\gl(V)$ be the endomorphism Lie algebra
spanned by the standard generators $\CE_{ij},\,1\leq i,j\leq\ell+1$
subjected to the
  relations:
 \be
  [\CE_{ij},\,\CE_{kl}]\,=\,\delta_{jk}\CE_{il}-\delta_{il}\CE_{kj}\,.
 \ee
Parabolic subalgebra $\frak{p}$ in $\gl_{\ell+1}$ is defined as a
subalgebra satisfying $\fb_-\subseteq\frak{p}\subset\gl_{\ell+1}$,
where $\fb_-\subset\gl_{\ell+1}$ is the Borel subalgebra
 spanned by $\{\CE_{ij},\,1\leq j\leq i\leq\ell+1\}$.
 Let $\fn_+\subset\gl_{\ell+1}$ be the nilpotent
subalgebra generated by $\{\CE_{ij},\,1\leq i<j\leq\ell+1\}$, so
that the following decomposition holds:
 \be\label{triang}
  \gl_{\ell+1}=\fb_-\oplus\fn_+\,=\,\fn_-\oplus\fh\oplus\fn_+\,,
 \ee
where $\fh\subset\fb$ is the Cartan subalgebra spanned by
$\CE_{nn},\,1\leq n\leq\ell+1$. Then
the Borel subalgebra $\fb_-\subset\gl_{\ell+1}$ is the minimal
parabolic subalgebra.  On the other hand the  maximal parabolic subalgebra
$\frak{p}_{1,\ell+1}\subset\gl_{\ell+1}$ is
generated by $\fb_-$ and
 \be
  \CE_{i,\,i+1}\,,\quad1\leq i<\ell\,.
  \ee

\subsection{Minimal parabolic $\gl_{\ell+1}$-Whittaker function}

The standard notion of the $\gl_{\ell+1}$-Whittaker function defined
for a reductive Lie algebra $\gl_{\ell+1}$ is a special case a more
general Whittaker function associated with a pair
$\mathfrak{p}\subset \gl_{\ell+1}$, for a parabolic subalgebra
$\mathfrak{p}$. Precisely, the standard Whittaker function thus
correspond to the case of minimal parabolic subalgebra (i.e. the
Borel subalgebra $\fb_-\subset\gl_{\ell+1}$).

For $\la=(\la_1,\ldots,\la_{\ell+1})\in\IC^{\ell+1}$, let
$(\rho_{\la},\,\CW_{\la})$ be the principal series representation of
the universal enveloping algebra $\CU(\gl_{\ell+1})$. Namely, for a
generic character $\chi_{\la}$ of the Borel subgroup $B_-\subset
GL_{\ell+1},\,\Lie(B_-)=\fb_-$ and the modular character
$\delta_{B_-}$, let $\CW_{\la}={\rm
Ind}_{B_-}^{GL_{\ell+1}}(\chi_{\la}\otimes\delta_{B_-}^{-1/2})$ the
induced representation. We consider the infinitesimal form ${\rm
Ind}_{\CU(\fb_-)}^{\gl_{\ell+1}}(\chi_{\la}\otimes\delta_{B_-}^{-1/2})$
of $\CW_{\la}$, and we choose the appropriately defined dual
$\CU(\gl_{\ell+1})$-module $(\rho^{\vee}_{\la},\CW^{\vee}_{\la})$.
Given the triangular decomposition \eqref{triang}, define after
\cite{Ko1},\cite{Ko2} the pair of Whittaker vectors
$\psi_R\in\CW_{\la}$ and $\psi_L\in\CW^{\vee}_{\la}$ to be the
generic characters of $\fn_-,\,\fn_+$:
 \be
  \rho^{\vee}_{\la}(\CE_{i+1,\,i})\psi_L\,=\,-\psi_L\,,\quad
  \rho_{\la}(\CE_{i,\,i+1})\psi_R\,=\,-\psi_R\,,\qquad1\leq
  i\leq\ell\,.
 \ee
We assume that the action of the subalgebra $\fh\subset\gl_{\ell+1}$
is integrated to the action of the maximal torus in
$GL_{\ell+1}(\IR)$. We also imply the existence of
$\gl_{\ell+1}$-invariant non-degenerate pairing $\<\,,\,\>$ between
the submodules $\CU(\gl_{\ell+1})\psi_L\subseteq\CW^{\vee}_{\la}$
and $\CU(\gl_{\ell+1})\psi_R\subseteq\CW_{\la}$. Then the
$\gl_{\ell+1}$-Whittaker function is given by the following matrix
element:
 \be
  \Psi_{\la_1,\ldots,\la_{\ell+1}}^{\gl_{\ell+1}}(e^{x_1},\ldots, e^{x_{\ell+1}})\,
  =\,e^{\sum\limits_{k=1}^{\ell+1}(\frac{\ell}{2}+1-k)x_k}\bigl\<\psi_L,\,
  \rho_{\la}\Big(e^{\sum\limits_{j=1}^{\ell+1}x_j\,\CE_{jj}}\Big)\psi_R\bigr\>\,.
  \ee
This matrix element may be expressed as the following
$\ell(\ell+1)$-fold integral \cite{GLO2} for ${\rm
Re}(\la_i)>0,\,1\leq i\leq\ell+1$:
 \be\label{GLnGiv}
  \Psi^{\gl_{\ell+1}}_{\la_1,\ldots,\la_{\ell+1}}(e^{x_1},\ldots,e^{x_{\ell+1}})
  =\,e^{x_{\ell+1}\sum\limits_{i=1}^{\ell+1}\la_i}\hspace{-2mm}
  \int\limits_{\IR_+^{2d}}\!
  \prod_{1\leq i\leq k\leq\ell}
  \frac{da_{k,i}}{a_{k,i}}\frac{db_{k,i}}{b_{k,i}}\,\,e^{-a_{k,i}-b_{k,i}}\\
  \times
  \prod_{1\leq i\leq k\leq\ell}(a_{k,i})^{\la_{k+1}}
  \prod_{k=1}^{\ell}(b_{k,k})^{\sum\limits_{i=1}^k\la_i}\\
  \times\prod_{i=1}^{\ell}\delta\Big(e^{-x_i+x_{i+1}}a_{\ell,i}b_{\ell,i}\,-\,1\Big)
  \prod_{1\leq i\leq k<\ell}\!\!
  \delta\Big(a_{k,i}b_{k,i}a_{k+1,\,i+1}^{-1}b_{k+1,\,i}^{-1}\,-\,1\Big)\,,
 \ee
where $d=\frac{\ell(\ell+1)}{2}$. For $\la_1=\ldots=\la_{\ell+1}=0$,
the integral \eqref{GLnGiv} originally appeared in \cite{Giv}. The
structure of the integrand in \eqref{GLnGiv}  might be elucidated by
invoking the type $A_{\ell}$ Gelfand-Zetlin graph:
 \be\label{GZgraph}
 \xymatrix{
  x_1 &&&&\\
  (\ell,1)\ar[u]^{a_{\ell,1}} & x_2\ar[l]_{b_{\ell,1}} &&&\\
  \vdots\ar[u] & \ddots\ar[u] & \ddots &&\\
%  (2,1)\ar[u]^{a_{21}} & \ddots\ar[l]_{b_{21}}\ar[u] & \ddots &
%  x_{\ell}\ar[l]_{b_{\ell,\ell-1}}\\
  (1,1)\ar[u]^{a_{11}} & %(2,2)\ar[l]^{b_{11}}\ar[u]^{a_{22}} &
  \ldots\ar[l]^{b_{11}} & (\ell,\ell)\ar[l]\ar[u]^{a_{\ell\ell}} &
  x_{\ell+1}\ar[l]^{b_{\ell\ell}}  }
 \ee
Namely, we attach  the integration variables $a_{k,i},\,b_{k,i}$ to
the edges, and add the arguments $x_i,\,1\leq i\leq\ell+1$ at the
boundary vertexes of the graph. The delta-factors in the integrand
of \eqref{GLnGiv} correspond to fixing products of variables
$a_{k,i}$ and $b_{k,i}$ along various paths on the diagram
\eqref{GZgraph}. Namely, we consider the $\ell$ paths starting
horizontally at $x_{i+1}$ then turning upwards to the adjacent
vertex $x_i$, and the $\ell(\ell-1)/2$ elementary box-shaped paths.

On the other hand, the integral \eqref{GLnGiv} can be identified
with the GG-integral of the form \eqref{INT1} associated with
$|I|=2d=\ell(\ell+1)$ (corresponding to the arrows in
\eqref{GZgraph}) and $|J|=d=\ell(\ell+1)/2$ (corresponding to the
vertices in \eqref{GZgraph}). Let us introduce the integration
variables $t_i\in\IR_+,\,i\in I$:
 \be
  \bigl\{t_i\,,\,\,1\leq i\leq\ell(\ell+1)\bigr\}\,
  =\,\{a_{k,i}\,,\,\,b_{k,i}\,:\quad1\leq i\leq k\leq\ell\}\,,
 \ee
and the two sets of arguments:
 \be\label{GLnGKZvar}
  \bigl\{y_{k,i}\,,\,\,z_{k,i}\,:\quad1\leq i\leq
  k\leq\ell\bigr\}\,\subset\,\IR^{2d}\,,\\
  \bigl\{\g_{k,i}\,,\,\,\nu_{k,i}\,:\quad
  1\leq i\leq k\leq\ell\bigr\}\,\subset\,\IC^{2d}\,,\qquad{\rm
  Re}(\g_{k,i})>0,\quad{\rm
  Re}(\nu_{k,i})>0\,.
 \ee
Consider the following integral of GG-type:
 \be\label{GLnGKZ}
  \Phi^{\gl_{\ell+1}}_{\g,\,\nu}(e^y,\,e^z)\,
  =\!\int\limits_{\IR_+^{2d}}\!\prod_{1\leq k\leq i\leq\ell}\!\!\Big\{
  \frac{da_{k,i}}{a_{k,i}}\frac{db_{k,i}}{b_{k,i}}\,\,
  a_{k,i}^{\g_{k,i}}b_{k,i}^{\nu_{k,i}}\,
  e^{-\,a_{k,i}-b_{k,i}}\\
  \times\prod_{i=1}^{\ell}
  \delta\Big(e^{-y_{\ell,i}-z_{\ell,i}}a_{\ell,i}b_{\ell,i}\,-\,1\Big)\\
  \times\prod_{1\leq i\leq k<\ell}\!\!
  \delta\Big(e^{-y_{k,i}-z_{k,i}+y_{k+1,\,i+1}+z_{k+1,\,i}}
  a_{k,i}b_{k,i}a_{k+1,\,i+1}^{-1}b_{k+1,\,i}^{-1}\,-\,1\Big)\,,
 \ee
associated with the following GKZ-data. Introduce the standard
orthonormal basis in $\IZ^{2d}$:
 \be
  \bigl\{\e^{k,i},\,\tilde{\e}^{\,k,i}\,:\quad
  1\leq i\leq k\leq\ell\bigr\}\,\subset\,\IZ^{2d}\,.
 \ee
the relation lattice,
 \be
  \IL=\Span\{l^{k,i}\,,\,\,1\leq i\leq k\leq\ell\}\subseteq\IZ^{2d}\,,
 \ee
is generated by
 \be\label{GLnL}
  l^{k,i}\,
  =\,\e^{k,i}+\tilde{\e}^{\,k,i}\,-\,\e^{k+1,\,i+1}-\tilde{\e}^{\,k+1,\,i}\,
  \in\,\IL\,,
 \ee
where $\e^{\ell+1,\,i+1}=\tilde{\e}^{\,\ell+1,\,i}=0$ is assumed.
Let $M\in\Mat_{d\times2d}$ be the relation matrix with rows given by
the generators $l^{k,i},\,1\leq i\leq k\leq\ell$. The defining
matrix, $A\in\Mat_{d\times2d}(\IZ)$, should satisfy the
orthogonality relation \eqref{ORTH}, $AM^{\top}=0\in\Mat_{d\times
d}(\IZ)$. One might choose the rows of $A$ to be the following:
 \be\label{GLnA}
  \a^{k,1}\,
  =\,\sum_{j=1}^k(\e^{k,j}-\tilde{\e}^{\,k,j})\,,\qquad1\leq k\leq\ell\,,\\
  \a^{k,i}\,
  =\,\sum_{j=1}^{i-1}\e^{k-j,\,i-j}\,
  +\,\sum_{j=i}^k(\e^{k,j}-\tilde{\e}^{\,k,j})\,,\qquad1<i\leq k\leq\ell\,.
 \ee
Let us stress that GG-hypergeometric function \eqref{GLnGKZ}
effectively depends only on the $d$ independent variables,
 \be\label{DeltaVar}
  R_i(y,z)\,:=\,y_{\ell,i}+z_{\ell,i}\,,\qquad1\leq i\leq\ell\,,\\
  B_{k,i}(y,z)\,:=\,y_{k,i}+z_{k,i}-y_{k+1,\,i+1}-z_{k+1,\,i}\,,
  \qquad  1\leq i\leq k<\ell\,,
 \ee
which is obvious from the integral representation \eqref{GLnGKZ}.
This may be attributed to the set of linear equations \eqref{TORUS1}
satisfied by the integral \eqref{GLnGKZ}.

Now let us introduce the following restriction of the
GG-hypergeometric function \eqref{GLnGKZ},
 \be\label{GLnGKZres}
  {}^{\rm res}\Phi^{\gl_{\ell+1}}_{\g,\,\nu}(e^x)\,:
  =\,\Phi^{\gl_{\ell+1}}_{\g,\,\nu}(e^y,\,e^z)\Big|_{\CL}\,,
 \ee
onto the $\ell$-dimensional {linear subspace $\CL\subset\IR^d$ with
the coordinates $x_i-x_{i+1},\,1\leq i\leq\ell$. Precisely, the
subspace $\CL$ is} defined by the following constraints on the
independent coordinates \eqref{DeltaVar} in $\IR^d$:
 \be\label{GLnGKZvarID}
  R_i(y,z)\,=\,x_i-x_{i+1}\,,\qquad1\leq i\leq\ell\,,\\
  B_{k,i}(y,z)\,=\,0\,,
  \qquad  1\leq i\leq k<\ell\,.
 \ee
\begin{lem} For the following special values of spectral variables from
\eqref{GLnGKZvar},
  \be\label{GLnGKZres2}
   \g_{k,i}(\lambda)=\la_{k+1}\,,\quad1\leq i\leq k\leq\ell\,;\qquad
   \nu_{k,k}(\lambda)=\sum_{i=1}^k\la_i\,,\quad1\leq k\leq\ell\,,\\
  \nu_{k,i}(\lambda)=0,\quad 1\leq i<k\leq\ell\,,
 \ee
the restricted GG-hypergeometric function \eqref{GLnGKZres} is
expressed through the $\gl_{\ell+1}$-Whittaker function
\eqref{GLnGiv}  as follows:
 \be\label{GLnGKZid}
 {}^{\rm res} \Phi^{\gl_{\ell+1}}_{\g(\lambda),\,\nu(\lambda)}(e^x)\,=
  e^{-\sum\limits_{i=1}^{\ell+1}\la_i\,
    x_{\ell+1}}\Psi^{\gl_{\ell+1}}_{\la}(e^x)\,.
 \ee
\end{lem}
\proof The identification directly follows by comparing
explicit integral representations \eqref{GLnGiv} and \eqref{GLnGKZ}. $\Box$

The difference in the exponential pre-factor may be taken into
  account  by considering a slightly extended GKZ data and taking a
  limit of the resulting hypergeometric function.
 Let us extend the sets $I,J$ defined by the graph \eqref{GZgraph} by
adding one element to each  $J$ and $I$, and introduce the extended
relation and defining  matrices
$\wh{M},\,\wh{A}\in\Mat_{(d+1)\times(2d+1)}(\IZ)$:
 \be\label{GLrel}
  \wh{J}:=J\sqcup\{d+1\}\,,\qquad\wh{I}:=I\sqcup\{2d+1\}\,,\qquad
  d\,=\,\frac{\ell(\ell+1)}{2}\,,\\
  \wh{M}:
  =\Big(\begin{smallmatrix}
   M&&0\\&&\\0&&1
  \end{smallmatrix}\Big)\,,\quad
  \wh{A}:
  =\Big(\begin{smallmatrix}
   A&&0\\&&\\0&&0
  \end{smallmatrix}\Big)\,,\qquad
  \wh{A}\wh{M}^{\top}=0\in\Mat_{(2d+1)\times(2d+1)}(\IZ)\,,
 \ee
so that the extended relation lattice $\wh{\IL}$ is spanned by its
rows $l^{k,i},\,(k,i)\in J$ and $l^{d+1}=(0,\ldots,0,1)$. Then the
matrix element \eqref{wavefun} associated with the GKZ data
\eqref{GLrel} for $\g_*\in\IC,\,{\rm Re}(\g_*)>0$ reads
 \be\label{GLnGKZ1}
  \wh{\Phi}_{\g,\nu;\,\g_*}(e^y,e^z;\,e^{y_*})\,
  =\!\int\limits_{\IR_+^{2d+1}}\!\!\!\frac{dt_{2d+1}}{t_{2d+1}}\,
  \delta(e^{-y_*}t_{2d+1}\,-\,1)\,t_{2d+1}^{\g_*}\,e^{-t_{2d+1}}\\
  \times\prod_{1\leq k\leq i\leq\ell}\!\!\Big\{
  \frac{da_{k,i}}{a_{k,i}}\frac{db_{k,i}}{b_{k,i}}\,\,
  a_{k,i}^{\g_{k,i}}b_{k,i}^{\nu_{k,i}}\,
  e^{-\,a_{k,i}e^{y_{k,i}}-b_{k,i}e^{z_{k,i}}}\\
  \times\prod_{i=1}^{\ell}
  \delta\Big(e^{-y_{\ell,i}-z_{\ell,i}}a_{\ell,i}b_{\ell,i}\,-\,1\Big)\\
  \times\prod_{1\leq i\leq k<\ell}\!\!
  \delta\Big(e^{-y_{k,i}-z_{k,i}+y_{k+1,\,i+1}+z_{k+1,\,i}}
  a_{k,i}b_{k,i}a_{k+1,\,i+1}^{-1}b_{k+1,\,i}^{-1}\,-\,1\Big)\\
  =\,e^{\g_*y_*\,-\,e^{y_*}}\times
  \Phi^{\gl_{\ell+1}}_{\g,\,\nu}(e^y,e^z)\,.
 \ee
Therefore, the $\gl_{\ell+1}$-Whittaker function \eqref{GLnGiv} can
be obtained as the following limit of the GG-hypergeometric function
\eqref{GLnGKZ1}:
 \be
  \Psi^{\gl_{\ell+1}}_{\la_1,\ldots,\la_{\ell+1}}(e^x)\,
  =\,\lim_{\e\to-\infty}\Big(
  {}^{\rm res}\wh{\Phi}_{\g(\lambda),\nu(\lambda);\,\e^{-1}\g_*}(e^x,\,e^{\e
  y_*})\Big)
  \Big|_{y_*=x_{\ell+1}\atop \g_*=\la_1+\ldots+\la_{\ell+1}}\,,
 \ee
 where ${}^{\rm res}\wh{\Phi}$ is obtained by imposing the
 restrictions \eqref{GLnGKZvarID}. Thus $\gl_{\ell+1}$-Whittaker
 function belongs to a partial compactification of the space of
 GG-hypergeometric functions.

\subsection{Maximal  parabolic $\gl_{\ell+1}$-Whittaker function}

Recall from \cite{GLO4} the definition of generalized Whittaker
function associated with maximal
 parabolic subalgebra $\frak{p}_{1,\ell+1}\subset \gl_{\ell+1}$
 (called the $(1;\ell+1)$-Whittaker function
 in \cite{GLO4}).  To the maximal parabolic subalgebra
 $\frak{p}_{(1,\ell+1)}$ we associate the  following decomposition of $\gl_{\ell+1}$
 \be\label{MAXtriang}
  \gl_{\ell+1}\,
  =\,\fh^{(1;\,\ell+1)}\,\oplus\,\fn_-^{(1;\,\ell+1)}\oplus\fn_+\,,
  \ee
where  $\fn^{(1;\,\ell+1)}_-\subset\fb_{-}$ is the
$\ell(\ell+1)/2$-dimensional subalgebra generated by
 \be\label{MAXn}
  \CE_{\ell+1,\,1}\,;\qquad
  \CE_{m+1,\,m}\,,\quad\CE_{m+1,\,m+1}\,,\quad1\leq m<\ell\,.
  \ee
and $\mathfrak{h}^{(1,\ell+1)}$ is the $(\ell+1)$-dimensional
commutative subalgebra spanned by
 \be\label{MAXh}
  h_1^{(1;\,\ell+1)}=\CE_{11}\,,\qquad
  h_2^{(1;\,\ell+1)}=\CE_{22}+\ldots+\CE_{\ell+1,\,\ell+1}\,,\\
  h_{m+1}^{(1;\,\ell+1)}=\CE_{\ell+1,\,m}\,,\quad1<m\leq\ell\,.
  \ee

For $\la=(\la_1,\ldots,\la_{\ell+1})\in\IC^{\ell+1}$, let
$(\rho_{\la},\,\CW_{\la})$ be the principal series representation of
$\CU(\gl_{\ell+1})$. Then we choose the appropriately defined dual
module $(\rho^{\vee}_{\la},\CW^{\vee}_{\la})$. Given the
decomposition \eqref{MAXtriang}, we introduce the pair of vectors
 $\psi_L^{(1;\,\ell+1)}\in\CW^{\vee}_{\la}$ and $\psi_R\in\CW_{\la}$ to be
the generic characters of the subalgebras $\fn^{(1;\,\ell+1)}_-$ and
$\fn_+$, where $\fn^{(1;\,\ell+1)}_-$ is generated by \eqref{MAXn}
(see \cite{GLO4},\cite{O}):
 \be\label{PellPsiL}
  \rho^{\vee}_{\la}(\CE_{\ell+1,\,1})\psi_L^{(1;\,\ell+1)}\,
  =\,-(-1)^{\frac{\ell(\ell-1)}{2}}\,\psi_L^{(1;\,\ell+1)}\,,\\
  \rho^{\vee}_{\la}(\CE_{21})\psi_L^{(1;\,\ell+1)}
  =\rho^{\vee}_{\la}(\CE_{m+1,\,m})\psi_L^{(1;\,\ell+1)}=0\,,\quad1<m<\ell\,;
 \ee
and
 \be\label{PellPsiR}
  \rho_{\la}(\CE_{i,\,i+1})\psi_R\,=\,-\psi_R\,,\qquad1\leq
  i\leq\ell\,.
 \ee
We imply the existence of $\gl_{\ell+1}$-invariant non-degenerate
pairing $\<\,,\,\>$ between the submodules
$\CU(\gl_{\ell+1})\psi_L^{(1;\,\ell+1)}\subseteq\CW^{\vee}_{\la}$
and $\CU(\gl_{\ell+1})\psi_R\subseteq\CW_{\la}$. Then the maximal
parabolic Whittaker function is defined by the following matrix
element:
 \be\label{WhittMax}
  \Psi_{\la}^{(1;\,\ell+1)}(e^{x_1},e^{x_2},\ldots , e^{x_{\ell+1}})\,
  =\,e^{x_1-x_{\ell+1}}\<\psi_L^{(1;\,\ell+1)},\,
  \rho_{\la}\Big(e^{\sum_j \,x_jh^{(1;\,\ell+1)}_j}\Big)\psi_R\>\,.
 \ee
The expression \eqref{WhittMax} implies that the action of the part
of subalgebra $\fb_-\subset\gl_{\ell+1}$ may be integrated to the
action of the corresponding group. In the following we consider the
restricted maximal parabolic Whittaker function
\cite{GLO4},\cite{O}, for ${\rm Re}(\la_i)>0,\,1\leq i\leq\ell+1$
given by
 \be\label{PellMatel}
  \Psi_{\la}^{(1;\,\ell+1)}(e^x)\,:=\,\Psi_{\la}^{(1;\,\ell+1)}(e^x,1\,\ldots,1)\,
  =\,e^x\<\psi_L^{(1;\,\ell+1)},\,\rho_{\la}(e^{x\,\CE_{11}})\psi_R\>\\
  =e^{\la_{\ell+1}x}\!\int\limits_{\IR_+^{\ell}}\!\prod_{i=1}^{\ell}
  \frac{dt_i}{t_i}\,t_i^{\la_i-\la_{\ell+1}}\,e^{-t_i}\,
  e^{-e^x\prod\limits_{i=1}^{\ell}t_i^{-1}}\\
  =\,\!\!\!\int\limits_{\IR_+^{\ell+1}}\!\!\prod_{i=1}^{\ell+1}
  \frac{dt_i}{t_i}\,t_i^{\la_i}\,e^{-t_i}\,
  \delta\Big(e^{-x}\prod\limits_{i=1}^{\ell+1}t_i\,-1\,\Big)\,.
 \ee
Therefore, we only shall assume that the action of the generator
$h_1^{(1;\,\ell+1)}=\CE_{11}\in\fh^{(1;\,\ell+1)}$ integrates to the
action of the corresponding one-parameter subgroup in
$GL_{\ell+1}(\IR)$.

To put the restricted maximal parabolic Whittaker function
\eqref{PellMatel} into the framework  of GG-hyper-~geometric
functions, consider the following instance of GKZ-system with
$|J|=1$, $|I|=N=\ell+1,\,m=\ell$ and the following matrices
$A\in\Mat_{\ell\times(\ell+1)}(\IZ)$ and
$M\in\Mat_{1\times(\ell+1)}(\IZ)$
 \be
  M=(1\ldots1)\,,\qquad
  A
  =\left(\begin{smallmatrix}
  1 & -1 & 0 & \ldots & 0\\
  0 & \ddots & \ddots & \ddots & \vdots\\
  \vdots & \ddots & \ddots & \ddots & 0\\
  0 & \ldots & 0 & 1 & -1
  \end{smallmatrix}\right),\\
  A^{\top}=(\a_1,\ldots,\a_{\ell})\,,\qquad
  \a_j=\e_j-\e_{j+1}\in\IZ^{\ell+1},\qquad
  MA^{\top}=0.
 \ee
Here $\{\e_j:\,j\in I\}\subset\IZ^{\ell+1}$ is the standard basis.
Let $\fL_{\ell+1}$ be the Lie algebra \eqref{LA}, and for
$\la=(\la_1,\ldots,\la_{\ell+1})\in\IC^{\ell+1},\,{\rm
Re}(\la_i)>0$, let $\pi_{\la}$ be the $\CU(\fL_{\ell+1})$-module
$\CV_{\la}$ \eqref{rep} modeled on $\CS(\IR_+^{\ell+1})$. Vector
$\phi_R\in\CV_{\la}$ and covector $\phi_L\in\CV^{\vee}_{\la}$ read
from \eqref{psiLsol}:
 \be
  \phi_R(t)\,=\,\prod_{i\in I}e^{-t_i}\,,\qquad
  \<\phi_L,\,\varphi\>\,
  =\,\!\!\int\limits_{\IR_+^{\ell+1}}\!\!\!
  \prod_{i\in I}\frac{dt_i}{t_i}\,
  t_i^{\la_i}\,\delta\Big(\prod_{j\in I}t_j\,-\,1\Big)\,
  \varphi(t)\,,
 \ee
for an arbitrary test function $\varphi(t)\in\CS(\IR_+^{\ell+1})$.
Then the matrix element \eqref{wavefun} in the
$\fL_{\ell+1}$-representation $\CV_{\la}$  affords the following
integral representation:
 \be\label{Pmatel}
 \Phi^{(1;\,\ell+1)}_{\la_1,\ldots,\la_{\ell+1}}(e^{y_1},\ldots,e^{y_{\ell+1}})\,
  =\,\bigl\<\phi_L,\,\pi_{\la}\Big(e^{\sum_iy_i\CH_i}\Big)\phi_R\bigr\>
 \ee
 \be\label{Pell}
  =\!\!\!\int\limits_{\IR_+^{\ell+1}}\!\!
  \prod_{i\in I}\frac{dt_i}{t_i}\,t_i^{\la_i}\,
  e^{\la_iy_i\,-\,t_ie^{y_i}}\,
  \delta\Big(\prod_{j\in I}t_j\,-\,1\Big)\\
  =\!\!\!\int\limits_{\IR_+^{\ell+1}}\!\!
  \prod_{i\in I}\frac{dt_i}{t_i}\,t_i^{\la_i}\,e^{-\,t_i}\,
  \delta\Big(\prod_{j\in I}e^{-y_j}t_j\,-\,1\Big)\,.
 \ee
Let us introduce variables $u_i:=e^{y_i},\,i\in I$ and write
\eqref{Pmatel} as follows:
 \be\label{PellGKZ}
  \Phi^{(1;\,\ell+1)}_{\la_1,\ldots,\la_{\ell+1}}(u_1,\ldots,u_{\ell+1})\,
  =\!\!\!\int\limits_{\IR_+^{\ell+1}}\!\!
  \prod_{i\in I}\frac{dt_i}{t_i}\,t_i^{\la_i}\,e^{-\,t_i}\,
  \delta\Big(\prod_{j\in I}u_j^{-1}t_j\,-\,1\Big)\,.
 \ee
By the results of Section 2, (or via a  simple direct
  computation),  the GG-hypergeometric function
$\Phi^{(1;\,\ell+1)}_{\la}(u)$ satisfies the following instances of
equations \eqref{DEF1}-\eqref{TORUS1d}:
 \be\label{PeqGKZ}
  \prod_{i\in I}\Big(-\frac{\pr}{\pr u_i}+\frac{\la_i}{u_i}\Big)
  \cdot\Phi^{(1;\,\ell+1)}_{\la}(u)\,
  =\,\Phi^{(1;\,\ell+1)}_{\la}(u)\,;\\
  \prod_{i\in I}e^{\pr_{\la_i}}\cdot
  \Phi^{(1;\,\ell+1)}_{\la}(u)\,
  =\,\prod_{i\in I}u_i\,\Phi^{(1;\,\ell+1)}_{\la}(u)\,;\\
  \Big(u_i\frac{\pr}{\pr u_i}-u_j\frac{\pr}{\pr u_j}\Big)\cdot
  \Phi^{(1;\,\ell+1)}_{\la}(u)\,=\,0\,,\qquad i\neq j\,;\\
  \Big\{e^{\pr_{\la_i}}\,-\,e^{\pr_{\la_j}}\Big\}\cdot
  \Phi^{(1;\,\ell+1)}_{\la}(u)\,=\,0\,,\qquad i\neq j\,,
 \ee
and the following instance of \eqref{DDeq}:
 \be\label{PellDD}
  \Big\{-u_i\frac{\pr}{\pr u_i}\,+\,\la_i\Big\}\cdot
  \Phi^{(1;\,\ell+1)}_{\la}(u)\,
  =\,e^{\pr_{\la_i}}\cdot\Phi^{(1;\,\ell+1)}_{\la}(u)\,.
 \ee

The first line of  equations in \eqref{PeqGKZ} can be written in the
following form:
 \be\label{Peq}
  \Big\{\prod_{i\in I}\Big(-u_i\frac{\pr}{\pr u_i}+\la_i\Big)\,
  -\,\prod_{i\in I}u_i\Big\}\cdot
  \Phi^{(1;\,\ell+1)}_{\la}(u)\,=\,0\,,
  \ee
As a consequence of the third line of equations in \eqref{PeqGKZ},
the function $\Phi^{(1;\,\ell+1)}_{\la}(u)$
  reduces to  a function in one variable
 \be\label{PellRed}
  \Phi^{(1;\,\ell+1)}_{\la_1,\ldots,\la_{\ell+1}}
   (e^{y_1},\ldots,e^{y_{\ell+1}})=
  \widetilde{\Phi}^{(1;\,\ell+1)}_{\la_1,\ldots,\la_{\ell+1}}(e^{x})\,,\qquad
  x=y_1+\ldots+y_{\ell+1},
 \ee
satisfying the reduced form of \eqref{Peq}:
 \be\label{PeqRed}
  \Big\{\prod_{i\in I}\Big(-\frac{\pr}{\pr x}+\la_i\Big)\,-\,e^x\Big\}\cdot
  \widetilde{\Phi}^{(1;\,\ell+1)}_{\la}(e^x)\,=\,0\,.
 \ee
Indeed, the equation \eqref{Peq} is a simple consequence of the
first line in \eqref{PeqGKZ}. By the third equation in
\eqref{PeqGKZ}, the function \eqref{Pell} depends on a single
variable $x=y_1+\ldots+y_{\ell+1}$ :
 \be
  \Phi^{(1;\,\ell+1)}_{\la_1,\ldots,\la_{\ell+1}}
   (e^{y_1},\ldots,e^{y_{\ell+1}})\,
   =\!\!\int\limits_{\IR_+^{\ell+1}}\!\!
  \prod_{i\in I}\frac{dt_i}{t_i}\,t_i^{\la_i}\,e^{-\,t_i}\,
  \delta\Big(e^{-\sum\limits_{j=1}^{\ell+1}y_j}\prod_{j\in
    I}t_j\,-\,1\Big)\,.
  \ee
Thus we may introduce the function
$\widetilde{\Phi}_{\la_1,\ldots,\la_{\ell+1}}(e^x)$, which by
\eqref{PellRed} has the following integral representation:
 \be\label{PellR}
  \widetilde{\Phi}^{(1;\,\ell+1)}_{\la_1,\ldots,\la_{\ell+1}}(e^x)\,:=\!\!\int\limits_{\IR_+^{\ell+1}}\!\!
  \prod_{i\in I}\frac{dt_i}{t_i}\,t_i^{\la_i}\,e^{-\,t_i}\,
  \delta\Big(e^{-x}\prod_{j\in I}t_j\,-\,1\Big)\,,
  \ee
and satisfies the equation \eqref{PeqRed} by construction.

Finally, we may identify the maximal parabolic Whittaker function $
\Psi^{(1;\,\ell+1)}_{\la_1,\ldots,\la_{\ell+1}}(e^x)$ defined by
\eqref{PellMatel} with the reduced GG-hypergeometric function
\eqref{PellRed} as follows:
 \be
  \Psi^{(1;\,\ell+1)}_{\la_1,\ldots,\la_{\ell+1}}(e^x)\,
   =\,\widetilde{\Phi}^{(1;\,\ell+1)}_{\la_1,\ldots,\la_{\ell+1}}(e^{x})\,.
 \ee

\noindent {\small {\bf A.A.G.} {\sl Laboratory for Quantum Field
Theory
and Information},\\
\hphantom{xxxx} {\sl Institute for Information
Transmission Problems, RAS, 127994, Moscow, Russia};\\
\hphantom{xxxx} {\it E-mail address}: {\tt anton.a.gerasimov@gmail.com}}\\
\noindent{\small {\bf D.R.L.} {\sl Laboratory for Quantum Field
Theory
and Information},\\
\hphantom{xxxx}  {\sl Institute for Information
Transmission Problems, RAS, 127994, Moscow, Russia};\\
\hphantom{xxxx} {\sl Moscow Center for Continuous Mathematical
Education,\\
\hphantom{xxxx} 119002,  Bol. Vlasyevsky per. 11, Moscow, Russia};\\
\hphantom{xxxx} {\it E-mail address}: {\tt lebedev.dm@gmail.com}}\\
\noindent{\small {\bf S.V.O.} {\sl
 School of Mathematical Sciences, University of Nottingham\,,\\
\hphantom{xxxx} University Park, NG7\, 2RD, Nottingham, United Kingdom};\\
\hphantom{xxxx} {\it E-mail address}: {\tt oblezin@gmail.com}}


\begin{thebibliography}{15}


\bibitem[BCFKvS]{BCFKvS} V. Batyrev, I. Ciocan-Fontanine, B. Kim,
  D. van Straten,   {\it Mirror
symmetry and toric degenerations of partial flag manifolds}, Acta
Math. 184 (2000) 1--39; Preprint {\tt arXiv:9803108}.


\bibitem[GG]{GG} I.M. Gelfand, M.I. Graev,
  {\it GG-functions and their relation with general hypergeometric functions},
  Rus. Math. Surveys 52:4 (1997) 3--48.

\bibitem[GGZ]{GGZ} I.M. Gelfand, M.I. Graev, A.V. Zelevinsky,
  {\it Holonomic systems of equations and series of hypergeometric type},
  Sov. Math. Dokl. 36 (1988) 5--10.

\bibitem[GKZ]{GKZ} I.M. Gelfand, M.M. Kapranov, A.V. Zelevinsky,
  {\it Hypergeometric functions and toric varieties},
  Funct. Anal. Appl. 23:2 (1989) 12--26.


\bibitem[G]{G} A. Gerasimov, {\it Archimedean Langlands duality
and exactly solvable quantum systems}, in Proc. ICM Seoul 2014, Vol.
3, 1097--1121.


\bibitem[GLO1]{GLO1} A.~Gerasimov, D.~Lebedev, S.~Oblezin Baxter operator
and Archimedean Hecke algebra, Commun. Math. Phys. 284 (2008)
867--896; Preprint {\tt arXiv:0706.3476}.

\bibitem[GLO2]{GLO2} A.~Gerasimov, D.~Lebedev, S.~Oblezin,
New integral representations of Whittaker functions for classical
groups, Rus. Math. Surveys 67:1 (2012) 3--96; Preprint {\tt
arXiv:0705.2886}.

\bibitem[GLO3]{GLO3} A.~Gerasimov, D.~Lebedev, S.~Oblezin, Archimedean
L-factors and topological field theories I,II, Commun. Number Theory
and Physics 5 (2011) 57--101, 101--134; Preprints {\tt
arXiv:0906.1065}, {\tt arXiv:0909.2016}.

\bibitem[GLO4]{GLO4} A. Gerasimov, D. Lebedev, S. Oblezin,
{\it Parabolic Whittaker functions and topological field theories
I}, Commun. Number Theory Phys. 5:2 (2011) 135--201; Preprint {\tt
arXiv:1002.2622}.

\bibitem[GLO5]{GLO5} A. Gerasimov, D. Lebedev, S. Oblezin,
{\it On a matrix element representation of special functions
associated with toric varieties}, in Nankai Symposium on
Mathematical Dialogues. Celebrating the 110th anniversary of the
birth of Prof. S.-S. Chern, Springer, 2022; Preprint {\tt
arXiv:2112.15013}.


\bibitem[Giv]{Giv} A.~Givental, {\it Stationary Phase Integrals,
Quantum Toda Lattices, Flag Manifolds and the Mirror Conjecture}.
Topics in Singularity Theory, AMS Transl. Ser. 2, 180, AMS,
Providence RI 1997, 103--115; Preprint {\tt arXiv:9612001}.


\bibitem[Ha]{Ha} M.~Hashizume,
{\it Whittaker functions on semi-simple Lie groups}, Hiroshima
Math.J., 12, (1982), 259--293.

\bibitem[Ja]{Ja}H.~ Jacquet, {\it Fonctions de Whittaker associ\'ees aux
groupes de Chevalley}, Bull.Soc.Math. France, 95 (1967), 243--309.

\bibitem[Ko1]{Ko1} B. Kostant, {\it On Whittaker vectors and representation theory},
Invent. Math. 48:2 (1978) 101--184.

\bibitem[Ko2]{Ko2} B. Kostant, {\it Quantization and representation theory},
in Proc. of Symposium on Representations and Lie groups, Oxford
1977, London Math. Soc. Lect. Notes Series, 1979, 34, 287--316.

\bibitem[O]{O} S.~Oblezin, On parabolic Whittaker functions II,
Cent. Eur. J. Math. 10:2 (2012) 543--558; Preprint {\tt
arXiv:1107.2998}.

\bibitem[OP]{OP} M.A. Olshanetsky, A.M. Perelomov, {\it Quantum integrable
systems related to Lie algebras}, Phys. Rep. 94:6 (1983), 313--404.


\bibitem[Sch]{SCH}G. Schiffmann, {\it Int\'egrales d'entrelacement et
fonctions de Whittaker}, Bull. Soc. Math. France 99 (1971), 3--72.


\end{thebibliography}
\end{document}